\documentclass[12pt]{amsart}

\usepackage[margin=3cm]{geometry}

\usepackage{amssymb,amsthm,amsmath}
\usepackage{hyperref}
\usepackage[alphabetic]{amsrefs}
\usepackage{dsfont,enumerate}

\newtheorem{theorem}{Theorem}[section]
\newtheorem{corollary}[theorem]{Corollary}
\newtheorem{lemma}[theorem]{Lemma}
\newtheorem{proposition}[theorem]{Proposition}

\theoremstyle{remark}
\newtheorem{remark}[theorem]{Remark}

\numberwithin{equation}{section}

\newcommand{\N}{\mathbb{N}}

\newcommand{\C}{\mathbb{C}}
\newcommand{\Crm}{\mathrm{C}}
\newcommand{\LSC}{\mathrm{Lsc}}
\newcommand{\Cu}{\mathrm{Cu}}

\newcommand{\K}{\mathcal{K}}

\newcommand{\waysubset}{\subseteq\!\subseteq}
\DeclareMathOperator{\rank}{rank}

\title{The Cuntz semigroup of some spaces of dimension at most two}

\author{Leonel Robert}

\address{Department of Mathematics\\ 
University of Louisiana at Lafayette\\
Lafayette, USA.}
\email{lrobert@louisiana.edu}

\begin{document}

\begin{abstract}
It is shown that the Cuntz semigroup  of a space with dimension at most two, and with second cohomology of its compact subsets equal to zero, is isomorphic to the ordered semigroup of lower semicontinuous functions on the space with values in the natural numbers with the infinity adjoined. 
This computation is then used to obtain the Cuntz semigroup of all compact surfaces. 
A converse to the first computation is also proven: if the Cuntz semigroup of a separable C*-algebra is isomorphic
to the lower semicontinuous functions on a topological space with values in the extended natural numbers, then the C*-algebra is commutative up to stability, and its spectrum satisfies the dimensional and cohomological conditions mentioned above.


\end{abstract}


\maketitle

\section{Introduction}
The Cuntz semigroup has become a useful tool in the general study of C*-algebras and in particular in their classification
(see  \cite{cowardelliottivanescu},  \cite{nccw}, \cite{toms}). It is thus necessary to compute examples and find ways to describe it. The Cuntz semigroup of commutative
C*-algebras is of particular interest, since they form the building blocks of various general classes of C*-algebras such as  AH and  ASH C*-algebras. 
In this paper we give a concrete description of the Cuntz semigroup of two classes of spaces: 2-dimensional spaces  with vanishing second cohomology class, and compact surfaces. Although these computations can be deduced from \cite[Proposition 5.4]{robert-tikuisis} (which covers all spaces of dimension at most 3) the proofs given here are more direct and the resulting computations are stated more explicitly. 
  
Let $X$ be a locally compact Hausdorff space of covering dimension at most 2. 
Let us denote by $\Cu(X)$ the Cuntz semigroup of the C*-algebra $\Crm_0(X)$. 
Let $\overline\N$ denote the set of natural numbers with 0 and $\infty$ adjoined. Let  
$\LSC_\sigma(X,\overline \N)$ denote the ordered semigroup of lower semicontinuous functions 
$f\colon X\to \overline\N$ such that the open sets $f^{-1}((n,\infty])$ are $\sigma$-compact for all $n\geq 0$ (we will drop the subscript $\sigma$ if this condition is redundant or not required).

\begin{theorem}\label{trivial}
Let $X$ be as above. Suppose that $\check H^2(K)=0$ (the \v Cech cohomology with integer coefficients) for every $K\subseteq X$
compact. Then the ordered semigroup $Cu(X)$ is isomorphic to $Lsc_\sigma(X,\overline \N)$.
 \end{theorem}

Theorem \ref{trivial} covers all spaces of dimension 0 or 1.
The case $\dim(X)=0$ was obtained previously by Perera in \cite{perera} 
and the case $X=[0,1]$ is implicit in Ivanescu's proof of \cite[Theorem 2.3]{ivanescu}. 
As pointed out above, Theorem \ref{trivial} can be deduced 
with a little work from \cite[Proposition 5.4]{robert-tikuisis}.

Let us now state the computation of the Cuntz semigroup of a compact surface.
Let $X$ be a compact surface.
Let $\mathrm V(X)$ denote the semigroup of isomorphism classes of finitely generated projective modules over $\Crm(X)$.
Let $(\LSC(X,\overline\N)\sqcup \mathrm V(X))/\!\!\sim$ be the quotient of the disjoint union of
$\LSC(X,\overline\N)$ and $\mathrm V(X)$ obtained by identifying the class $n[\Crm(X)]$ in $\mathrm V(X)$ 
with the constant function
$n$ in $\LSC(X,\overline\N)$. Let us define on this set 
an order and an addition operation. Inside $\LSC(X,\overline\N)$ and $\mathrm V(X)$, we retain the order and addition with which these sets are endowed. 
Let $f\in \LSC_\sigma(X,\overline\N)$ be a non-constant
function and let $[P]\in \mathrm V(X)$. Then
 $f+[P]$  is defined as the function
$f+\rank([P])\in \LSC(X,\overline\N)$. Let the order relation be  $f\leq [P]$ if $f\leq \rank([P])$ and $[P]\leq f$ if $\rank([P])\leq f$ (otherwise $f$ and $[P]$ are not comparable).

\begin{theorem}\label{sphere}
Let $X$ be a compact surface. Then  
\begin{enumerate}[(i)]
\item
$Cu(X)$ is isomorphic to the ordered semigroup $(Lsc(X,\overline 
\N)\sqcup V(X))/\!\!\sim$.
\item
The relations of Cuntz equivalence and of isomorphism agree for the countably generated Hilbert C*-modules over $C(X)$.
\end{enumerate}
\end{theorem}

Part (ii) of the previous theorem should be compared with \cite[Example 1]{robert-tikuisis}, where a 2-dimensional CW complex is exhibited for which the relations of Cuntz equivalence and of isomorphism do not agree. (On the other hand, \cite[Example 1]{robert-tikuisis} is covered by Theorem \ref{trivial}.) 
En route to proving Theorems \ref{trivial} and \ref{sphere}, we establish some results  of interest in themselves. 
Proposition \ref{iso} below belongs to the general theory of Hilbert C*-modules and seems to not have been pointed out before: if $I$ is closed two-sided ideal of a C*-algebra and $G$ and $H$ are countably generated (right) Hilbert C*-modules over that algebra such that $G/GI\cong H/HI$, then
$H\oplus GI\cong G\oplus HI$. This improves on \cite[Theorem 1.2]{crs}, where under the same assumptions and 
with $I$ $\sigma$-unital it is shown that $H\oplus \ell_2(I)\cong G\oplus \ell_2(I)$. 

The last section of the paper is devoted to proving the following converse to Theorem \ref{trivial}:
\begin{theorem}\label{converse}
Let $A$ be a separable C*-algebra. Suppose that $A\cong Lsc(X,\overline{\N})$ for some locally compact Hausdorff space $X$. Then $A$ is stably isomorphic to $C_0(X)$. Furthermore, $X$ has dimension at most 2 and $\check H^2(K)=0$ for all $K\subseteq X$, with $K$ compact.
\end{theorem}

\section{Results on Hilbert C*-modules}
Let $A$ be a C*-algebra and let $H$ be  a Hilbert C*-module over $A$ (all modules are assumed to be right modules). Given  a closed two-sided ideal $I$ of $A$, 
$HI$  denotes the closed submodule  $\{h\cdot i\mid h\in H,i\in I\}$. The quotient $H/HI$ is assumed to be endowed with its natural structure of $A/I$-Hilbert C*-module. (We will some times simply say Hilbert module rather than Hilbert C*-module, always meaning the same thing.) 

If $A=\Crm_0(X)$,  we will denote by $H_x$  the Hilbert space $H/H\Crm_0(X\backslash\{x\})$ (referred to as  the fibre of $H$ at $x$). We will denote by $\ell_2(\Crm_0(X))$  the $\Crm_0(X)$-Hilbert module of infinite sequences $(f_i)$, with $f_i\in \Crm_0(X)$ such that $\sum_{i=1}^\infty f_i^*f_i$ is convergent in norm. 
Finally, $\mathcal K$ will denote the C*-algebra of compact operators over a separable Hilbert space.

\begin{lemma}\label{likeKR}
Let $A$ be a C*-algebra. Let $a,b\in A_+$ be of norm at most 1 and let $\epsilon$ be such that $\|a-b\|<\epsilon$.
Then there is $y$ such that $y^*y=(a-\epsilon)_+$, $yy^*\leq b$, and $\|yy^*-b\|<3\epsilon+2\sqrt{\epsilon}$.
\end{lemma}
\begin{proof}
The proof of this lemma follows closely the first part of the proof
of \cite{kirchberg-rordam}*{Lemma 2.2}.
Let $\epsilon_1$ be such that $\|a-b\|<\epsilon_1<\epsilon$.
Let $e\in C^*(a)$ be a positive contraction such that $e(a-\epsilon_1)e=(a-\epsilon)_+$. Set $b^{1/2}e=x$ and let $x=v|x|$
be the polar decomposition of $x$, where $v\in A^{**}$. From $a-\epsilon_1\leq b$ we get $(a-\epsilon)_+\leq ebe=x^*x$,
and so the element $v(a-\epsilon)_+^{1/2}$ belongs to $A$. Set $v(a-\epsilon)_+^{1/2}=y$. We have $y^*y=(a-\epsilon)_+$ and
\[
yy^*=v(a-\epsilon)_+v^*\leq vx^*xv^*=xx^*=b^{1/2}e^2b^{1/2}\leq b.
\]
Also, 
\[\|yy^*-b\|\leq \|v(a-\epsilon)_+ - vx^*xv^*\|+\|b^{1/2}e^2b^{1/2}-b\|.\]
With a few of simple computations it may be shown that the right-hand side is bounded
by $3\epsilon+2\sqrt{\epsilon}$. This is done
using that $e,v,a$ and $b$ are contractions, and that $\|a^{1/2}-b^{1/2}\|<\sqrt{\epsilon}$
(this is implied by $\|a-b\|<\epsilon$).
\end{proof}

For each  $f\in (\Crm_0(X)\otimes \mathcal K)_+$ and $x\in X$, let us  denote by  $\sigma_1(f)(x),\sigma_2(f)(x),\ldots$ the eigenvalues of $f(x)$ arranged in decreasing order.

\begin{proposition}
Let $X$ be a locally compact Hausdorff space of dimension at most 2.
Let $B$ be a $\sigma$-unital hereditary subalgebra of $C_0(X)\otimes \mathcal K$.
Then the set of strictly positive elements $f\in B_+$ such that all the non-zero eigenvalues
of $f(x)$ have multiplicity 1 for all $x\in X$ is a dense subset of $B_+$. 
\end{proposition}

\begin{proof}
The proof relies on a Baire category argument.
Let $g\in B_+$ be a strictly positive element. For each $n=1,2,\dots$, let $U_n=\{x\in X\mid \mathrm{rank}(g(x))\geq n\}$.
Observe that the sets $U_n$ depend only on the algebra $B$ and not on the choice of the
strictly positive element $g$ (since $\rank g(x)=\sup \{\rank h(x)\mid h\in B_+\}$).

Let $K\subseteq U_n$ be compact and $\epsilon>0$. 
Let us denote by $V(n,K,\delta)$ the subset of $B_+$ of elements $f\in B_+$ 
such that 
\begin{enumerate}[(1)]
\item
$\sigma_{n}(f)(x)>\sigma_{n+1}(f)(x)$ for all $x\in K$,

\item
$\|g-fb\|<\delta$ for some $b\in B$.
\end{enumerate}

Let us show that $V(n,K,\delta)$ is open and dense in $B_+$. Since neither
of the properties (1) and (2) is destroyed by a sufficiently small perturbation of $f$,  $V(n,K,\delta)$ is an open subset of $B_+$. Let us show that any $h\in B_+$ may be approximated 
by elements of $V(n,K,\delta)$. Since the strictly positive elements form a  dense in $B_+$, 
we may assume that $h$ is strictly positive.  In particular,  $\sigma_n(h)(x)$ is non-zero for all $x\in K\subseteq U_n$. Let 
\[
\epsilon_0=\min_{x\in K}\sigma_n(h)(x).\] 
Let $\epsilon>0$
and assume that $\epsilon<\epsilon_0/2$. By \cite[Theorem 5.3]{benoit} (see \cite[Theorem 1]{choi-elliott} for the metrizable case) we may approximate $h$ by a positive element $h'$ in $\Crm_0(X)\otimes \mathcal K$ such that the $n$-th largest eigenvalue of $h'(x)$ has multiplicity one  and is non-zero for all $x\in X$.  
Let us choose such an $h'$ such that $\|h-h'\|<\epsilon$  and 
\begin{align}\label{ntheigen}
|\sigma_n(h)(x)-\sigma_n(h')(x)| <\frac {\epsilon_0}{2}
\end{align} 
for all $x\in K$.
By Lemma \ref{likeKR}, there
exists $v\in \Crm_0(X)\otimes \mathcal K$ such that $(h'-\epsilon)_+=v^*v$ and $vv^*\leq h$. Since $\epsilon<\epsilon_0/2$, \eqref{ntheigen} implies that the $n$-th largest eigenvalue of $ (h'-\epsilon)_+(x)$ is simple and non-zero for all $x\in K$. Thus, the same is true for  
 $n$-th largest  eigenvalue of $vv^*$. 
Let us choose $\epsilon_1>0$ small enough such that the $n$-th largest eigenvalue of 
$\epsilon_1 g+vv^*$ has multiplicity one
for all $x\in K$. Set $\epsilon_1 g+ vv^*=f$. The function $f$ satisfies
that $\epsilon_1 g\leq f\leq \epsilon_1g+h$, whence $f$ is a strictly positive element of $B$. In particular,  $f$ satisfies (2) for any value of $\delta>0$. Finally,  by letting $\epsilon$ and $\epsilon_1$ be arbitrarily small, $f$ will get arbitrarily close to $h$ (observe that $vv^*$ is close to $h$ by Lemma \ref{likeKR}). Hence $V(n,K,\delta)$ is dense in $B_+$.

For each $n$  the set  $U_n$ is $\sigma$-compact. Thus,  there are compact sets $(K_{n,m})_{m=1}^\infty$ such that $\bigcup_{m=1}^\infty K_{n,m}=U_n$. By the Baire category theorem, $\bigcap_{n,m,l=1}^\infty V(n,K_{n,m},\frac{1}{2^{l}})$ is a dense subset of $B_+$. It is easily verified that the elements of this intersection have the properties stated in the theorem. 
\end{proof}

\begin{theorem}\label{cormaind2}
Let $X$ be a locally compact Hausdorff space of dimension at most 2. Let $H$ be a countably
generated Hilbert C*-module over $C_0(X)$. Then
\[
H\cong \bigoplus_{n=1}^\infty P_nC_0(U_n),
\]
where $(U_n)_{n=1}^\infty$ is a decreasing sequence of $\sigma$-compact open subsets of $X$, and
$P_n$ are finitely generated projective modules over $C_b(U_n)$ of constant dimension 1 over their fibres.
\end{theorem}
\begin{proof}
Since $H$ is countably generated there is a positive element $f\in \Crm_0(X)\otimes \mathcal K$ such that $H\cong \overline{f\ell_2(C_0(X))}$. By the previous theorem applied to the hereditary subalgebra generated by $f$, we may assume that the non-zero  eigenvalues of $f(x)$ have multiplicity 1  for all $x\in X$. For each $n=1,2,\dots$, let $U_n=\{x\in X\mid \sigma_n(f)(x)>0\}$. Then $f$ has the form $f=\sum_{n=1}^\infty p_n \sigma_n(f)$, where $p_n\in C_b(U_n)\otimes \mathcal K$ is such that $p_n(x)$ is the spectral projection of $f(x)$ associated to the eigenvalue $\sigma_n(f)(x)$ for all $x\in U_n$. Since the non-zero eigenvalues of $f(x)$ have multiplicity 1, $p_n(x)$ is a rank 1 projection for all $n$ and all $x\in U_n$. Also,  the elements $(p_n\sigma_n(f))_{i=1}^\infty$ (which are defined on all $X$) are mutually orthogonal.
We have
\[
\overline{f\ell_2(C_0(X))}\cong \bigoplus_{n=1}^\infty \overline{p_n \sigma_n(f)\ell_2(C_0(X))}.
\]
So, setting $p_n\cdot \ell_2(C_b(U_n))=P_n$ we get the desired result.
\end{proof}

The second result used in the computation of the Cuntz semigroup applies
to countably generated Hilbert C*-modules over an arbitrary C*-algebra. 

\begin{proposition}\label{iso}
Let $I$ be a closed two-sided ideal of a C*-algebra $A$. Let $H$ and $E$
be countably generated Hilbert C*-modules over $A$ such that $H/HI$ is 
isomorphic to $E/EI$ as $A/I$-Hilbert C*-modules, and let $\phi\colon H/HI\to E/EI$ 
be an isomorphism between them. Then the Hilbert modules $H\oplus EI$ and $E\oplus HI$ are
isomorphic and there is an isomorphism $\Phi\colon H\oplus EI\to E\oplus HI$ that 
lifts $\phi$.
\end{proposition}

\begin{proof}
By the nonconmutative Tietze extension Theorem (see \cite[Theorem 3]{crs}), there is $\psi\colon H\to E$
that lifts $\phi$, and $\|\psi\| \le 1$. Let $\Phi\colon H \oplus E \to E \oplus H$ be given by the matrix
\[
\Phi=
\begin{pmatrix}
\psi & \sqrt{1-\psi\psi^*}\\
\sqrt{1-\psi^*\psi} & \psi^*
\end{pmatrix}.
\]
It is easily verified that $\Phi$ is an isomorphism of Hilbert C*-modules. Since $\psi$ lifts
an isomorphism from $H/HI$ to $E/EI$, we have that the image of $1-\psi\psi^*$ is
contained in $EI$ and the image of $1-\psi^*\psi$ is contained in $HI$. The same is true
for the square roots of these operators. This implies that the restriction of $\Phi$ to
$H \oplus EI$ is an isomorphism from $H \oplus EI$ to $E \oplus HI$.
\end{proof}

\begin{remark} It is shown in \cite[Theorem 1.1]{crs} that if $I$ is $\sigma$-unital, and $H$ and $E$ are countably generated Hilbert modules such that $E/EI$ is  Cuntz smaller than  $H/HI$, then $E\oplus HI$ is Cuntz smaller than $H\oplus EI$ (see the next section for the definition of Cuntz equivalence). This leads us to the following question: Suppose that $E/EI$ embeds isometrically in $H/HI$ (assume also that $I$ is $\sigma$-unital if necessary). Does $E\oplus HI$ embed in $H\oplus EI$? 
\end{remark}

\begin{corollary}
Let $I$ and $J$ be $\sigma$-unital closed two-sided ideals of $A$ and
let $H$ be a countably generated Hilbert C*-module over $A$. Then 
\[HI\oplus HJ\cong H(I+J)\oplus H(I\cap J).\]
\end{corollary}
\begin{proof}
This follows from Proposition \ref{iso} applied to the countably generated Hilbert modules
$H(I+J)$ and $HJ$, which are isomorphic after passing to the quotient
by $I$.
\end{proof}

\section{Proof of Theorems \ref{trivial} and \ref{sphere}}
Let us briefly recall the Hilbert modules picture of the Cuntz semigroup given by  Coward, Elliott, and Ivanescu. The reader is referred to their paper  \cite{cowardelliottivanescu} for the details of this construction. Let $A$ be a C*-algebra. Let $H$ be a Hilbert module over $A$ and let $F\subseteq H$ be a Hilbert submodule of $H$. Let us write $F\waysubset H$ if there exists $T\in \mathcal K(H)_+$ such that $T|_F=\mathrm{id}_F$. Now let  $G$ and $H$ be countably generated Hilbert modules. Let us say that $H$ is Cuntz smaller than $G$ if for every $F\waysubset H$ there exists $F'\waysubset G$
such that $F\cong F'$ (isomorphism of Hilbert modules). This relation is denoted by $H\preceq G$. Furthermore,  let us say that $G$ and $H$ are Cuntz equivalent, and denote this by $H\sim G$, if they are each Cuntz smaller than the other. Observe that isomorphic Hilbert modules are clearly Cuntz equivalent but the converse is known not to always hold. 
The Cuntz semigroup of the C*-algebra $A$ is then defined as the set of Cuntz equivalence classes of countably generated Hilbert modules over $A$. This set is denoted by $\Cu(A)$.
Addition in $\Cu(A)$ is defined by 
\[
[H]+[G]=[H\oplus G],
\] 
and the order is defined by $[H]\leq [G]$ if $H$ is Cuntz smaller than $G$.  

If $A=\mathrm C_0(X)$, we will abbreviate $\Cu(\mathrm C_0(X))$ to $\Cu(X)$.

The following lemma clarifies the role of the cohomological condition in Theorem
\ref{trivial}.

\begin{lemma}\label{clarifiesH2}
Let $X$ be a $\sigma$-compact, locally compact Hausdorff space. Suppose that $\check H^2(K)=0$ for any $K\subseteq X$ compact.
Let $P$ be a finitely generated projective $C_b(X)$-Hilbert module of constant dimension 1 on its fibres. 
Then $PC_0(X)\sim C_0(X)$.
\end{lemma} 

\begin{proof} (cf. \cite[p. 250]{robert-tikuisis}.)
Let $U$ be an open set compactly contained in $X$, i.e., such that $U\subseteq K\subseteq  X$ for some compact set $K$.
Consider the module $P|_K:=P/P\Crm_0(X\backslash K)$. The first Chern class (with values in $\check {\mathrm{H}}^2$) classifies line bundles up to isomorphism (see \cite[Theorem 16.3.4]{husemoller}). Since $\check{\mathrm H}^2(K)=0$, the line bundle associated to $P_K$ is trivial and so $P_K\cong \Crm(K)$.
Since $P\Crm_0(U)$ can be naturally identified, by the quotient map $P\rightarrow P|_K$, with $P_K\Crm_0(U)\cong \Crm_0(U)$, we conclude that $P\Crm_0(U)\cong \Crm_0(U)$. In particular, $[P\Crm_0(U)]=[\Crm_0(U)]$ in $\Cu(X)$. Passing to the supremum over an increasing sequence of compactly contained open subsets $U$ that cover $X$, we get that $[P\Crm_0(X)]=[\Crm_0(X)]$, as desired.
\end{proof}

\begin{proof}[Proof of Theorem \ref{trivial}]
Let $H$ be a countably generated Hilbert module over $\Crm_0(X)$.
By Theorem \ref{cormaind2}, $H$  is isomorphic -- whence Cuntz equivalent -- to a module of the form $\bigoplus_{i=1}^\infty P_i\Crm_0(U_i)$. Furthermore, applying the previous lemma to each $U_i$,  we get that $H$ is Cuntz equivalent to $\bigoplus_{i=1}^\infty \Crm_0(U_i)$. Here $(U_i)_{i=1}^\infty$ is a decreasing sequence of
$\sigma$-compact open subsets of $X$. 

Let us show that the map $\rank\colon \Cu(X)\to \LSC_\sigma(X)$, given by
\[
\rank([H])(x)=\dim H_x, \mbox{\quad for }x\in X,
\]
is an isomorphism. Indeed, we may assume that  $H=\bigoplus_{i=1}^\infty \Crm_0(U_i)$ for some decreasing sequence of $\sigma$-compact open sets $(U_i)_{i=1}^\infty$. Then $\rank([H])=\sum_{i=1}^\infty \mathds{1}_{U_i}$, where $\mathds{1}_{U_i}$ denotes the characteristic function of $U_i$. Since every function $f\in \LSC_\sigma(X,\overline N)$ can be represented in the form $\sum_{i=1}^\infty \mathds{1}_{U_i}$ (with $U_i=f^{-1}((i,\infty])$), the map $\rank(\cdot)$ is surjective. If the sequence of open sets $(U_i)_{i=1}^\infty$ is required to be decreasing then the representation of $f$ in this form is unique. Hence $\rank(\cdot)$ is also injective. It is easily verified that $\rank(\cdot)$ preserves order and addition. 
\end{proof}

The following lemma contains the topological facts about surfaces that we need in the proof of Theorem \ref{sphere}.
\begin{lemma}
Let $X$ be a compact surface and let
 $U$ be a proper, non-empty, open subset of $X$. 
\begin{enumerate}[(i)]
\item 
 Let $P$ be a rank 1 projective Hilbert $\Crm_b(U)$-module. Then $P\Crm_0(U)\cong \Crm_0(U)$.
 \item
 Let $Q$ be a rank 1 projective Hilbert $\Crm(X)$-module.  Then $Q\oplus \Crm_0(U)\cong \Crm(X)\oplus \Crm_0(U)$.
 \end{enumerate} 
\end{lemma}
\begin{proof}
The crucial properties of a compact surface that we need are 
\begin{enumerate}
\item
every proper open subset of a compact surface is homotopic to a 1-dimensional space,
\item
every proper cosed subset of a compact surface  has vanishing $\check{\mathrm H}^2$ group.
\end{enumerate}
Property (1) is proven -- in more general form -- in \cite{vrabec} (see also \cite[Proposition 1.1]{repovs} for the case $X=S^2$). Property (2) is well known and can be derived as follows:  on one hand $\check{\mathrm H}^2(X\backslash\{x\})=0$ for any $x\in X$ (by property (1)). On the other hand, $\check H^3(X\backslash\{x\},C)=0$ for any closed subset $C\subset X\backslash\{x\}$, since $X\backslash \{x\}$ is 2-dimensional. Thus, by the exact sequence in cohomology, 
$\check{\mathrm H}^2(C)=0$.
Both (1) and (2) imply that proper open and proper closed subsets of a compact surface can only have trivial complex line bundles, since the first Chern class -- with values in $\check{\mathrm H}^2$ -- classifies line bundles.  Thus, with $U$ and $P$ as in part (i) of the lemma, we must have $P\cong \Crm_b(U)$, and so $P\Crm_0(U)\cong \Crm_0(U)$. Also, the Hilbert modules $Q$  and $\Crm(X)$ of part(ii) become isomorphic  after passing to the quotient by $\Crm_0(U)$ (i.e., after restricting them to $X\backslash U$). Therefore, by Proposition \ref{iso},  
\[
Q\oplus \Crm_0(U)\cong \Crm(X)\oplus Q\Crm_0(U)\cong \Crm(X)\oplus \Crm_0(U).\qedhere
\] 
\end{proof}

\begin{proof}[Proof of Theorem \ref{sphere}]
Let $H$ be countably generated Hilbert C*-module $H$ over $\Crm(X)$. By Theorem \ref{cormaind2}, we may assume that 
$H=\bigoplus_{i=1}^\infty P_i\Crm_0(U_i)$.
Suppose that there are indices $i$ such that the subset $U_i$ is proper and nonempty and let $i_0$ be the smallest such
index. Then $P_j\Crm_0(U_j)\cong \Crm_0(U_j)$ for all $j\geq i_0$ by part (i) of the previous lemma. 
Furthermore, by (ii) of the previous  lemma, 
\[
P_{1}\oplus P_2\oplus\dots \oplus P_{i_0-1}\oplus \Crm_0(U_{i_0})\cong \Crm(X)\oplus \Crm(X)\oplus\dots \oplus \Crm(X)\oplus \Crm_0(U_{i_0}).
\]
Hence, $H\cong \bigoplus_{i=1}^\infty \Crm_0(U_i)$. Suppose on the other hand  that the subsets $U_i$ are never proper and nonempty.
If $U_i=X$ for all $i$ then $H\cong \ell_2(\Crm(X))$ (by \cite{dixmier-douady}*{Theorem 5}). If $U_i$ is eventually empty then $H$ is a finitely
generated projective module over $\Crm(X)$. Summarizing, every countably generated Hilbert C*-module over $\Crm(X)$ is isomorphic to either (1) a finitely generated projective module, (2) a Hilbert module of the form $\bigoplus_{i=1}^\infty \Crm_0(U_i)$, where $(U_i)_{i=1}^\infty$ is a decreasing sequence of open subsets of $X$. The assertion in part (ii) of the theorem now clearly follows from this. For if two $\Crm(X)$-Hilbert modules $H$ and $G$ are Cuntz equivalent, then $\rank(H)=\rank(G)$. 
In the case that they have constant finite rank, they are both in case (1). That is, they are Cuntz equivalent finitely generated projective modules. It is known that in this situation Cuntz equivalence and isomorphism agree (this follows from the definition of Cuntz equivalence and the fact that projections in $\Crm(X)\otimes \mathcal K$ are finite.
In case (2) the modules are again isomorphic, since the open sets $(U_i)_{i=1}^\infty$ can be  recovered from the rank function
(see the proof of Theorem \ref{trivial}).

Let us now compute $\Cu(X)$.
Consider the map $\Cu(X)\to \LSC_\sigma(X,\overline \N)\sqcup V(X)/\!\!\sim$ 
given by
\begin{align}\label{isomap}
\begin{array}{rcl}
[\bigoplus_{i=1}^\infty \Crm_0(U_i)] &\mapsto& \sum_{i=1}^\infty \mathds{1}_{U_i},\\
\left[Q\right] &\mapsto& [Q].
\end{array}
\end{align}
Let $Q$ be a finitely generated projective module of dimension $m$ over its fibres and let $H=\bigoplus_{i=1}^\infty \Crm_0(U_i)$, 
with at least one $U_i$ 
proper and nonempty; say $U_{i_0}$ is the first such set. Since the module $Q\oplus H$ has non-constant dimension function (equal to $m+\rank([H])$), we must have
$Q\oplus H\cong \Crm(X)^m\oplus \bigoplus_{i=1}^\infty \Crm_0(U_i)$. 
This shows that the map defined by \eqref{isomap} is linear.

Let $Q$ and $H$ be as before. 
Suppose that $\rank([Q])=m$ and $\rank([H])\geq m$. Notice that $i_0\geq m$. 
By Theorem \ref{cormaind2}, $Q\cong P_1\oplus P_2\dots \oplus P_m$,
for some projective modules $P_i$, $i=1,2,\dots,m$, of constant dimension 1.
For every such
module $P_i$, we have $P_i\oplus \Crm(U_{i_0})\cong \Crm(X)\oplus \Crm(U_{i_0})$, and so $P_i$ is isomorphic to a submodule of 
$\Crm(X)\oplus \Crm(U_{i_0})$.
Hence, $P_1\oplus P_2\dots \oplus P_m$ is isomorphic to a submodule of $\Crm(X)^m\oplus \Crm(U_{i_0})$.
Since $i_0\geq m$, we have $[Q]\leq [H]$.

Suppose on the other hand that $\rank([Q])=m$ and $\rank([H])\leq m$. Then 
$H=\Crm_0(U_1)\oplus \dots \oplus \Crm_0(U_k)$ for some $k\leq m$. The modules
$P_1\Crm_0(U_1)\oplus \dots \oplus P_k\Crm_0(U_k)$ and $H$ are isomorphic, since they
have the same (non-constant) dimension function. On the other hand, 
$P_1\Crm_0(U_1)\oplus \dots \oplus P_k\Crm_0(U_k)$ is clearly a submodule of $Q$.
It follows that $[H]\leq [Q]$.
This completes proving that the map defined by \eqref{isomap} is an isomorphism of ordered semigroups.
\end{proof}

\section{Proof of Theorem \ref{converse}}
In this section we will make use of the positive elements picture of $\Cu(A)$. In this case the elements
of $\Cu(A)$ are equivalence classes of positive elements of $A\otimes \K$. Given $a,b\in (A\otimes\K)_+$,
we say that $a$ is Cuntz smaller than $b$, and write $a\preceq b$, if there exist $d_n\in A\otimes \K$ such that $d_n^*bd_n\to a$. We say that $a$ is Cuntz equivalent to $b$, and denote this by $a\sim b$, if $a\preceq b$ and $b\preceq a$. The Cuntz semigroup $\Cu(A)$ is defined as the set of Cuntz equivalence classes of positive elements. We shall denote the class of $a\in A\otimes\K$ by $[a]$. The reader is referred to \cite{ara-perera-toms} for a presentation of this construction of the Cuntz semigroup.
It is shown in \cite{cowardelliottivanescu} that the positive elements picture and the Hilbert modules picture indeed give the same ordered semigroup. 
 
Let us first show that the ideals of $A$ are in bijective correspondence with the ideals of $\Cu(A)$. 
A subsemigroup  $S\subseteq \Cu(A)$ is called an ideal if it is hereditary (i.e., $t\leq s$ and $s\in S$ imply 
$t\in S$) and closed under the suprema of increasing sequences.
\begin{lemma}
Let $A$ be a C*-algebra. For each closed 2-sided ideal $I$ of $A$, let 
\[
S_I=\{[a]\mid a\in (I\otimes \mathcal K)_+\}.
\]
Then $I\mapsto S_I$ is an order isomorphism from the lattice of  closed two sided
ideals of  $A$ to the lattice of ideals of $\Cu(A)$.
\end{lemma}
\begin{proof}
The lattice of closed two-sided ideals is unchanged when passing from $A$ to $A\otimes\K$. Thus, let us rename
$A\otimes \K$ by $A$ (i.e., assume that $A$ is stable) and define $S_I=\{[a]\mid a\in I_+\}$. 

Let us first show that $S_I$ is indeed an ideal. 
Observe that from the definition of the Cuntz comparison relation, if $a\preceq b$ and $b\in I_+$ then $aI_+$. 
In particular, if $a\sim a'$ then $a\in I$ if and only if $a'\in I$.
It thus follows that $S_I$ is a hereditary set. Suppose that $[a]=\sup_n [a_n]$, where $([a_n])_{n=1}^\infty$ is an increasing
sequence with $a_n\in I$ for all $n$. Then for each $\epsilon>0$ there exists $n$ such that $[(a-\epsilon)_+]\leq [a_n]$. So $(a-\epsilon)_+\in I$ and passing to the limit we get that $a\in I$. Thus, $S_I$ is closed under sequential suprema. The proof that $S_I$ is closed under addition is left to the reader.

It is clear that the map $I\mapsto S_I$ is order preserving. Suppose that $S_I=S_J$. Then for every $a\in I_+$ we have that $a\sim a'$, for some $a'\in J$. Thus, $a\in J$. This readily implies that $I=J$. Finally, given an ideal $S$ of $\Cu(A)$ define $I=\mathrm{Ideal}(\{a\in A\mid [a]\in S\})$. We have that $S\subseteq S_I$. Suppose that $a\in I_+$. Then for every $\epsilon>0$
there exist $d_i,d_2,\cdots,d_n\in A$ and $b_1,b_2,\cdots,b_n\in A_+$ such that $[b_i]\in S$ for all $i$ and
$(a-\epsilon)_+=\sum_{i=1}^n d_i^*b_id_i$. Thus, $[(a-\epsilon)_+]\leq \sum_{i=1}^n [b_i]\in S$. Passing to the supremum over all $\epsilon=1/1,1/2,1/3,\dots$ and using the $S$ is closed under such suprema we get that $[a]\in S$. Thus, $S=S_I$. 
\end{proof}

\begin{proof}[Proof of Theorem \ref{converse}]
Let us assume without loss of generality that $A$ is stable. Let $a\in A\otimes \mathcal K$ be a positive element
such that $[a]\in \Cu(A)$ is mapped by the isomorphism to the constant function $1\in \LSC(X,\overline{\N})$.
Observe that there is no $x\in \Cu(A)$ such that 
$2x\leq [a]$ with $x\neq 0$. By Glimm's lemma (see \cite[Proposition 3.10]{robert-rordam}), this implies that all irreducible representations of  $\overline{aAa}$ are 1-dimensional. That is, $\overline{aAa}\cong \Crm_0(Y)$ for some locally compact Hausdorff space $Y$.
Since $[a]$ is a  full element (i.e., $\infty\cdot [a]=\infty$), we have that $\overline{aAa}$ is a full hereditary subalgebra. Thus, by Brown's theorem, $A\cong \Crm_0(Y)\otimes \mathcal K$. By the previous lemma, $Y=\mathrm{Prim}(A)\cong X$. Thus, $A\cong \Crm_0(X)\otimes \mathcal K$.

We have $X$ such that $\Cu(X)\cong \LSC(X,\overline{\N})$. Let us argue that in fact this isomorphism
is given by the map $\Cu(X)\stackrel{\rank}{\longrightarrow}\LSC(X,\overline{\N})$. The map $\rank (\cdot)$ is clearly surjective. On the other hand,
if $[f],[g]\in \Cu(X)$ are such that $\rank([f])=\rank([g])$, then $[f]$ and $[g]$ agree after passing to any maximal quotient of $\Cu(X)$. In $\LSC(X,\overline{\N})$, this property implies that both elements must be equal. Since $\Cu(X)\cong \LSC(X,\overline{\N})$, we get that $[f]=[g]$. It follows
that $\rank(\cdot)$ is an isomorphism.

Next, let us show that $X$ has dimension at most 2. Passing to a quotient of $\Crm_0(X)$, we may assume that
$X$ is compact. Suppose that $\dim X>2$. By \cite[Theorem 1.9.3]{engelking}, this means that there exists 
a closed subset $Z\subseteq X$ and a continuous map $p\colon Z\to S^2$, such that $p$ cannot be extended to a continuous 
map on all $X$. Identifying the 2-dimensional sphere $S^2$ with the rank 1 projections of $M_2(\C)$, we may view $p$ as a rank 1 projection  in $M_2(\Crm(Z))$. Let us lift this projection to a positive element $f\in M_2(\Crm(X))$.  Set $\lambda\cdot 1_2+f=g$, where $\lambda\in \Crm_0(X\backslash Z)$ is a strictly positive function on $X\backslash Z$. Then $\rank([g])(x)=1$ for $x\in Z$ and $\rank([g])(x)=2$ for $x\in X\backslash Z$.
But $[g]$ cannot agree with the obvious Cuntz semigroup element having the same rank, namely, the element $[h]$ with 
$h=
\Big(
\begin{smallmatrix}
1 & 0\\
0 &\lambda
\end{smallmatrix}\Big)$. For suppose that this is the case. Since
$[1]\leq [h]$, we get that $[1]\leq [g]$. Thus, there exists a partial isometry $v\in M_2(\Crm(X))$ such that
$1=v^*v$ and $vv^*\in \overline{bM_2(\Crm(X))b}$. We must have that $vv^*|_Z=p$, since the left and the right hand side
are rank 1 projections and $vv^*|_Z\leq p$. But this implies that $vv^*$ is an extension of $p$ to all of $X$, which contradicts
our choice of $p$.

Finally, let us show that $\check{\mathrm H}^2(K)=0$ for any $K\subseteq X$ compact. Again, passing to the quotient we may assume that $K=X$. If $\check{\mathrm H}^2(X)\neq 0$, then $X$ has a non-trivial line bundle. Thus, there exists a rank 1 projection $q\in M_2(\Crm(X))$
which is not trivial. So $[q]\neq [1]$ in $\Cu(X)$, although they both have the same rank. This contradicts that the map
$\rank(\cdot)$ is an isomorphism from $\Cu(X)$ to $\LSC(X,\overline{\N})$.
\end{proof}

\textbf{Acknowledgments.} I thank Alex Karassev for clarifying  that 
\v Cech cohomology was the right cohomology to use in the hypotheses of Theorem \ref{trivial}.
I thank Du\v{s}an Repov\v{s} for sharing the unpublished reference \cite{vrabec} with me.  
I thank Luis Santiago and Franklin Vera Pacheco for  beneficial 
discussions of this paper.

\end{document}